\documentclass[11pt]{amsart}
\usepackage{amstext,amssymb,amsmath,amsbsy,dsfont,tikz}

 \usepackage[left=2.65cm,right=2.65cm,top=3.2cm,bottom=3.2cm]{geometry}

\usepackage{amsmath,amssymb,latexsym,dsfont}
\usepackage[small]{caption}
\usepackage{graphicx,color,mathrsfs,tikz}
\usepackage{subfigure}
\usepackage{cite}
\usepackage[colorlinks=true,urlcolor=blue,
citecolor=red,linkcolor=blue,linktocpage,pdfpagelabels,
bookmarksnumbered,bookmarksopen]{hyperref}
\usepackage[italian,english]{babel}
\usepackage{units}
\usepackage{enumitem}

\usepackage{hyperref}
\usepackage{cleverref}

\usepackage{tikz}
\usetikzlibrary{intersections}

\usepackage{amscd}
\usepackage{amsfonts}
\usepackage{indentfirst}
\usepackage{verbatim}
\usepackage{amsmath}
\usepackage{amsthm}
\usepackage{enumerate}
\usepackage{graphicx}
\usepackage{color}
\usepackage[OT1]{fontenc}
\usepackage[latin1]{inputenc}
\usepackage[english]{babel}
\usepackage{amssymb,esint}

\newtheorem{theorem}{Theorem}[section]
\newtheorem{lemma}{Lemma}[section]

\newtheorem{remark}{Remark}[section]

\usepackage[normalem]{ ulem }
\usepackage{soul}

\usepackage{tikz}
\usetikzlibrary{intersections}

\numberwithin{equation}{section}

\newcommand{\tr}{^\mathsf{T}}

\newcommand{\dsp}{\displaystyle}

\newcommand{\cF}{{\mathcal F}}

\newcommand{\eps}{\varepsilon}

\newcommand{\mR}{\mathbb{R}}

\newcommand{\cB}{\mathcal B}

\newcommand{\B}{\mathrm{B}}

\newcommand{\D}{{\mathcal D}}

\newcommand{\hD}{ {\mathcal O}}

\title[Finite-time stabilization in optimal time of hyperbolic systems]{Finite-time stabilization in optimal time of homogeneous quasilinear  hyperbolic systems in one dimensional space}
\author{Jean-Michel Coron}
\address[Jean-Michel Coron]{Sorbonne Universit\'{e}, Universit\'{e} de Paris, CNRS, INRIA,
\newline \indent 	Laboratoire Jacques-Louis Lions, \'{e}quipe Cage, Paris, France.}
\email{coron@ann.jussieu.fr.}

\author{Hoai-Minh Nguyen}
\address[Hoai-Minh Nguyen]{Ecole Polytechnique F\'ed\'erale de Lausanne, EPFL, 
\newline \indent  CAMA, Station 8,  CH-1015 Lausanne, Switzerland.}
\email{hoai-minh.nguyen@epfl.ch}

\begin{document}

\maketitle

\begin{abstract}
We consider the finite-time stabilization  of  homogeneous quasilinear hyperbolic systems with one side controls and  with nonlinear boundary condition at the other side.  We present time-independent feedbacks leading to the finite-time stabilization in any time larger than  the optimal time for the null controllability of the linearized system if  the initial condition is sufficiently small.  One of the key technical points is to establish the local well-posedness of quasilinear hyperbolic systems with nonlinear, non-local boundary conditions.
\end{abstract}

\section{Introduction and statement of the main result}

Linear hyperbolic systems in one dimensional space are frequently used
in modeling of many systems such as traffic flow, heat exchangers,  and fluids in open channels. The
stability and boundary stabilization of these hyperbolic systems
have been studied intensively in the literature, see  e.g.   \cite{BastinCoron} and the references therein.
In this paper, we investigate the finite-time stabilization in optimal time  of the following homogeneous,  quasilinear,  hyperbolic system in one dimensional space
\begin{equation}\label{Sys-1}
\partial_t w (t, x) =  \Sigma \big(x, w(t, x) \big) \partial_x w (t, x)  \mbox{ for } (t, x)  \in [0, + \infty) \times (0, 1).
\end{equation}
Here $w = (w_1, \cdots, w_n)\tr: [0, + \infty) \times (0, 1) \to \mR^n$, $\Sigma(\cdot, \cdot)$ is an   $(n \times n)$ real matrix-valued function defined in $[0, 1] \times \mR^n$. We assume that  $\Sigma (\cdot, \cdot)$ has  $m \ge 1$  distinct positive eigenvalues and $k = n - m \ge 1$  distinct negative eigenvalues. As usual, see e.g. \cite{CoronNg19}, we assume that, maybe after a change of variables,   $\Sigma(x, y)$ for $x \in [0, 1]$ and   $y \in \mR^n$ is of the form
\begin{equation}\label{form-A}
\Sigma(x, y) = \mbox{diag} \Big(- \lambda_1(x, y), \cdots, - \lambda_{k}(x, y),  \lambda_{k+1}(x, y), \cdots,  \lambda_{k+m}(x, y) \Big),
\end{equation}
where
\begin{equation}\label{relation-lambda}
-\lambda_{1}(x, y) < \cdots <  - \lambda_{k} (x, y)< 0 < \lambda_{k+1}(x, y) < \cdots \lambda_{k+m}(x, y).
\end{equation}
Throughout the paper, we assume
\begin{equation}\label{cond-lambda}
\mbox{$\lambda_i$  is of class $C^2$ with respect to $x$ and $y$    for $1 \le i \le n = k + m$.}
\end{equation}
Denote
$$
\mbox{$w_- = (w_1, \cdots, w_k)\tr $ and $w_+ = (w_{k+1}, \cdots, w_{k+m})\tr$.}
$$
The following types of boundary conditions and controls are considered. The boundary condition at $x = 0$ is given by
\begin{equation}\label{bdry-w-0}
w_-(t, 0)  = \B \big(w_+(t,  0) \big) \mbox{ for } t \ge 0,
\end{equation}
for some
$$
\mbox{$\B \in \big( C^2(\mR^m) \big)^k$ \mbox{with} $\B(0) = 0$,}
$$
and the boundary control at $x = 1$ is
\begin{equation}\label{bdry-w-1}
w_+(t, 1) = (W_{k+1}, \cdots, W_{k+m})\tr (t) \mbox{ for } t \ge 0,
\end{equation}
where $W_{k +1}, \dots, W_{k + m}$ are controls. In this work, we thus consider non-linear boundary condition at $x = 0$.

Set
\begin{equation}\label{def-tau}
\tau_i =  \int_0^1  \frac{1}{\lambda_i(x, 0)} \, dx  \mbox{ for } 1 \le i \le n,
\end{equation}
and
\begin{equation}\label{def-Top}
T_{opt} := \left\{ \begin{array}{cl}  \dsp \max \big\{ \tau_1 + \tau_{m+1}, \dots, \tau_k + \tau_{m+k}, \tau_{k+1} \big\} & \mbox{ if } m \ge k, \\[6pt]
\dsp \max \big\{ \tau_{k+1-m} + \tau_{k+1},  \tau_{k+2-m} + \tau_{k+2},  \dots, \tau_{k} + \tau_{k+m} \big\} &  \mbox{ if } m < k.
\end{array} \right.
\end{equation}

The main result of this paper is  the following result whose proof is given in the next section.

\begin{theorem}  \label{thm1} Define
\begin{equation}
\cB: = \Big\{B \in \mR^{k \times m}; \mbox{ such that  \eqref{cond-B-1} holds for  $1 \le i \le \min\{m-1, k\}$} \Big\},
\end{equation}
where
\begin{multline}\label{cond-B-1}
\mbox{ the $i \times i$  matrix formed from the last $i$ columns and the last $i$ rows of $B$  is invertible.}
\end{multline}
Assume that $B = \nabla \B(0) \in {\mathcal B}$.  For any $T > T_{opt}$, there exist $\eps > 0$ and a time-independent feedback control for \eqref{Sys-1}, \eqref{bdry-w-0}, and \eqref{bdry-w-1}  such that if the compatibility conditions $($at $x =0$$)$  \eqref{compatibility-0} and \eqref{compatibility-1} below hold for $w(0, \cdot)$, 
\begin{equation}\label{finite-T}
\left(\| w(0, \cdot) \|_{C^1([0, 1])} < \eps\right)\Rightarrow \left(w(T, \cdot) =0\right). 
\end{equation}
\end{theorem}

\begin{remark}\rm 1. The feedbacks constructed also lead to the well-posedness of the Cauchy problem for the closed loop system (see Lemma~\ref{lem-QL}) and to the following property: for every $\eta>0$, there exists $\delta>0$ such that, if the compatibility conditions (at $x =0$)  \eqref{compatibility-0} and \eqref{compatibility-1} below hold for $w(0, \cdot)$,
\begin{equation}\label{stability-T}
\left(\| w(0, \cdot) \|_{C^1([0, 1])} < \delta\right)\Rightarrow \left(\| w(t, \cdot) \|_{C^1([0, 1])} < \eta, \, \forall t\in [0,T]\right);
\end{equation}
see the proof of Lemma~\ref{lem-QL}. Hence, by  \eqref{finite-T} and \eqref{stability-T}, $0\in \big((C^1([0,1])\big)^n$ is stable for the closed-loop system and $0\in \big(C^1([0,1])\big)^n$ is finite-time stable in time $T$. 2. The feedbacks constructed in this article use additional $4m$ state-variables (dynamics extensions) to avoid imposing compatibility conditions at $x=1$. In particular \eqref{finite-T} and \eqref{stability-T} are understood with these  additional $4m$ state-variables. 
\end{remark}

In what follows, we denote, for $x \in [0,1]$ and $y \in \mR^n$,
$$
\Sigma_-(x, y)  = \mbox{diag} \Big(- \lambda_1(x, y), \cdots, - \lambda_{k}(x, y) \Big) \mbox{ and } \Sigma_+(x, y)  = \mbox{diag} \Big( \lambda_{k+1}(x, y), \cdots,  \lambda_{n}(x, y) \Big).
$$
The compatibility conditions considered in \Cref{thm1} are:  \begin{equation}\label{compatibility-0}
w_-(0, 0) = \B  \big( w_+(0, 0) \big)
\end{equation}
and
\begin{equation}\label{compatibility-1}
\Sigma_- \big(0, w(0, 0)\big) \partial_x w_-(0, 0)=\nabla \B \big( w_+(0, 0) \big)  \Sigma_+ \big(0, w(0, 0)\big) \partial_x w_+(0, 0).
\end{equation}

Null-controllability of hyperbolic systems with one side controls have been studied at least  from the work of David Russell \cite{Russell78} even for inhomogeneous systems, i.e., instead of  \eqref{Sys-1}, one considers
$$
\partial_t w(t, x) = \Sigma \big(x,w(t. x) \big) \partial_x w(t, x) + C\big(x, w(t, x) \big),
$$
for some $C \in \big( L^\infty([0, 1] \times \mR^n) \big)^{n \times n}$ with $C(x, 0) =0$.
For linear systems, i.e., $\Sigma(x, \cdot)$ and $C(x, \cdot)$ are constant for $x \in [0, 1]$ and $\B$ is linear ($\B (\cdot) = B \cdot$ with $B = \nabla \B(0)$),  the null-controllability  was established in \cite[Section 3]{Russell78} for the time $\tau_k + \tau_{k+1}$. Using backstepping approach, feedback controls leading to finite-time stabilization in the same time  were then  initiated by Jean-Michel Coron et al.  in \cite{CVKB13} for $m=k = 1$ and later developed   in  \cite{AM16, CHO17} for the general case. The set $\cB$ was introduced in \cite{CoronNg19} and the null-controllability for the  linear systems with $B \in \cB$ was established for $T> T_{opt}$ in \cite{CoronNg19, CoronNg19-2} (see also \cite{Weck} for the case $C$ diagonal) via the backstepping approach.   A tutorial introduction of backstepping approach can be found in \cite{Krstic08}.  In the  quasilinear case with $m \ge k$ and with the linear boundary condition at $x=0$,   the null controllability for any time greater than $\tau_k + \tau_{k+1}$ was established for $m \ge k$ by Tatsien Li in  \cite[Theorem 3.2]{Li00} (see also \cite{LiRao02}).

This work is concerned about homogeneous quasilinear hyperbolic systems with controls on one side, and with nonlinear boundary conditions on the other side: \eqref{Sys-1}, \eqref{bdry-w-0}, and \eqref{bdry-w-1}.  When the boundary condition is linear, the null-controllability was obtained by Long Hu  \cite{LongHu}  for $m \ge k$ at  any time greater than $\max\{\tau_{k+1}, \tau_{k} + \tau_{m+1}\}$ if initial data are  sufficiently small. In the linear case \cite{CoronNg19},  for $B \in \cB$, we  obtained time-independent feedbacks for  the null controllability at the {\it optimal} time $T_{opt}$  and showed  the optimality of $T_{opt}$.  Related exact controllability results can be also found in  \cite{LongHu, CoronNg19,HO19}. In this work, for $\nabla \B(0) \in \cB$, we present time-independent feedbacks leading to finite-time stabilization of \eqref{Sys-1}, \eqref{bdry-w-0}, and \eqref{bdry-w-1} in any time $T> T_{opt}$ provided that the initial data are sufficiently small.  It is easy to see  that $\cB$ is an open subset of  the set of (real) $k\times m$ matrices,  and the Hausdorff dimension of its complement is $\min\{k, m-1 \}$.

The feedbacks for \eqref{Sys-1}, \eqref{bdry-w-0}, and \eqref{bdry-w-1} are nonlinear and  inspired from the ones in \cite{CoronNg19}. The construction is  more complicated due to quasilinear nature of the system. We add auxiliary dynamics to fulfill the compatibility conditions at $x=1$ since $C^1$-solutions are considered. One of the key technical points is to establish the local well-posedness of quasilinear hyperbolic systems with nonlinear, non-local boundary conditions, which is interesting in itself.

\section{Proof of the main result}

This section containing two subsections is devoted to the proof of \Cref{thm1}. In the first subsection, we establish the local well-posedness of quasilinear hyperbolic systems with nonlinear, non-local boundary conditions. This implies in particular  the well-posedness for  the feedback laws given in the proof of \Cref{thm1} associated with \eqref{Sys-1} and  \eqref{bdry-w-0}. The proof of \Cref{thm1} is given in the second subsection.

\subsection{Preliminaries} The main result  of this section is \Cref{lem-QL} on the  well-posedness  for quasi-linear hyperbolic systems related to \eqref{Sys-1} and  \eqref{bdry-w-0}. The assumptions made  are guided by  our  feedback controls used in \Cref{thm1}.  We first consider the semilinear system, with $T>0$,
\begin{equation}\label{lem-SL-sys}
\left\{\begin{array}{cl}
\partial_t u(t, x) = A(t, x) \partial_x u(t,x) + f\big(t, x, u(t, x) \big)& \mbox{ in } [0, T] \times [0, 1],  \\[6pt]
u_{-}(t, 0)  = g \big(t, u_{+}(t, 0) \big)  & \mbox{ for } t \in [0, T],\\[6pt]
u_{+}(t, 1) = h \big(t,  u(t, \cdot), u_0 \big)& \mbox{ for } t \in [0, T], \\[6pt]
u(0, \cdot) = u_{0}(\cdot) &  \mbox{ in } [0, 1],
\end{array}\right.
\end{equation}
for
$$
A(t, x) = \mbox{diag} \big(-\lambda_1(t, x), \cdots, - \lambda_m(t, x), \lambda_{m+1}(t, x), \cdots, \lambda_{m + k}(t, x) \big),
$$
where
$$
-\lambda_1(t, x) <  \cdots <  - \lambda_m(t, x)< 0 <  \lambda_{m+1}(t, x) <  \cdots <  \lambda_{m + k}(t, x),
$$
and for $f: [0, T] \times [0, 1] \times \mR^n \to \mR^n$, $g: [0,  T] \times \mR^{m} \to \mR^k$, and $h : [0, T] \times \big(C^1([0, 1]) \big)^n \times \big(C^1([0, 1]) \big)^n  \to \mR^m$.

\medskip
We have

\begin{lemma} \label{lem-SL} Assume that 
$A$ is  of class $C^1$, $f$, and  $g$ are of class $C^2$, 
\begin{equation}\label{lem-SL-cd-01}
h(t, \varphi, u_0) = h_1(t, u_0) + h_2(t, \varphi, u_0) \mbox{ with } h_1, \, h_2 \mbox{ are of class $C^1$}, 
\end{equation}
\begin{equation}\label{lem-SL-cd-02}
\lim_{\eta \to 0} \sup_{\| u_0 \|_{C^1([0, 1])} \le \eta} \sup_{t >0} \Big( | h_1(t, u_0) | + | \partial_t h_1(t, u_0) | \Big) = 0, 
\end{equation}
\begin{equation}\label{lem-SL-cd-1}
 \quad f(t, x, 0) = g(t, 0) = h_2(t, 0, \cdot) = 0, 
\end{equation}
and  the following conditions hold, for some $C>0$, $a \in [0, 1)$, $1 \le p < + \infty$, and $\eps_0>0$,
\begin{multline}\label{lem-SL-cd-2}
|h (t, \hat \varphi, u_0) -  h (t, \varphi, u_0)| + |\partial_t h (t, \hat \varphi, u_0) - \partial_t h (t, \varphi, u_0)|  \\[6pt]
\le C \Big( \| \big( \hat \varphi - \varphi, \hat \varphi' - \varphi' \big) \|_{C^0([0, a])} +  \| \big( \hat \varphi - \varphi,  \hat \varphi' - \varphi'  \big) \|_{L^p(0, 1)} \Big),
\end{multline}
 for all $\varphi,\,  \hat \varphi, \,  u_0 \in \big( C^1([0, 1]) \big)^n$ with $\max\big\{ \| \hat \varphi \|_{C^1( [0, 1])}, \| \varphi \|_{C^1( [0, 1])}, \| u_0 \|_{C^1( [0, 1])} \big\} < \eps_0$, and
\begin{multline}\label{lem-SL-cd-3}
\left|\frac{d}{dt}h(s, \hat v(t, \cdot), u_0) |_{s = t} - \frac{d}{dt}h(s, v(t, \cdot), u_0) |_{s = t} \right| \\[6pt]
\le C \Big( \| \big( (\hat v - v) (t, \cdot),  \partial_t (\hat v - v) (t, \cdot), \partial_x (\hat v - v) (t, \cdot) \big) \|_{C^0([0, a])} \\[6pt] +  \| \big( (\hat v - v) (t, \cdot),  \partial_t (\hat v - v) (t, \cdot), \partial_x (\hat v - v) (t, \cdot)  \big) \|_{L^p(0, 1)} \Big),
\end{multline}
 for all $\hat v, v \in \big(C^1([0, T] \times [0, 1])\big)^n$ and $u_0 \in \big(C^1([0, 1]) \big)^n$ with $ \max\big\{ \| \hat v \|_{C^1([0, T] \times [0, 1])}, \|  v \|_{C^1([0, T] \times [0, 1])} \big\} \\< \eps_0$ and $\| u_0\|_{C^1([0, 1])} < \eps_0$.
There exists $\eps > 0$ such that for $u_0 \in \big( C^1([0, 1])\big)^n$ satisfying the compatibility conditions $($see \eqref{lem-SL-comp-1}-\eqref{lem-SL-comp-3} below$)$ with $\|u_0\|_{C^1([0, 1])} < \eps$,  there  is  a unique solution $u \in \big(C^1([0, T] \times [0, 1]) \big)^n$ of \eqref{lem-SL-sys}.
\end{lemma}

We recall the following definition of compatibility conditions for \eqref{lem-SL-sys}: $u_0 \in \big( C^1([0, 1])\big)^n$ is said to satisfy the compatibility conditions if
\begin{equation}\label{lem-SL-comp-1}
u_{0, -}(0) = g\big(0, u_{0, +}(0) \big), \quad u_{0, +}(1) = h\big(0, u_0, u_0 \big),
\end{equation}
\begin{multline}\label{lem-SL-comp-2}
\Big( A(0, 0) u_0'(0) + f\big(0, 0, u_0(0) \big) \Big)_- \\[6pt]
 = \partial_t g(0, u_{0, +}(0) ) + \partial_{y_+} g (0, u_{0, +}(0)) \Big( A(0, 0) u_0'(0) + f\big(0, 0, u_0(0) \big) \Big)_+,
\end{multline}
\begin{multline}\label{lem-SL-comp-3}
\Big( A(0, 1) u_0'(1) + f\big(0, 1, u_0(1) \big) \Big)_+\\[6pt]
 =  \partial_t h(0, u_0, u_0)  + \partial_{y} h (0, u_0, u_0)   \Big( A(0, \cdot)  u_0'(\cdot) + f\big(0, \cdot, u_0(\cdot) \big) \Big).
\end{multline}
Here and in what follows, the partial derivatives are taken with respect to the notations $f(t, x, y)$, $g(t, y_+)$, and  $h(t, y, u_0)$.

\begin{remark}  \rm The conditions $a < 1$ and $p< + \infty$ are crucial in \Cref{lem-SL}.
\end{remark}

\begin{proof}[Proof of \Cref{lem-SL}]  Set, for $ u \in \big(C([0, T] \times [0, 1])\big)^n$,
\begin{equation}\label{0-norm}
\| u \|_0 : = \max_{1 \le i \le n} \mathop{\max}_{(t, x) \in [0, T] \times [0, 1]} | e^{-L_1 t - L_2 x} u_i(t, x)|
\end{equation}
and, for $ u \in \big(C^1([0, T] \times [0, 1])\big)^n$,
\begin{equation}\label{1-norm}
\| u \|_1 : = \max\Big\{\| u \|_0, \|\partial_t u \|_0,  \|\partial_x u \|_0 \Big\},
\end{equation}
where $L_1$ and $L_2$ are two  large,  positive constants determined later.

Set
\begin{multline*}
\hD_\eps : = \Big\{ v \in \big(C^1([0, T] \times [0, 1])\big)^n \mbox{ with }  v(0, \cdot) = u_0,  \\[6pt]
\partial_t v(0, 1) = A(0, 1) u_0'(1) + f\big(0, 1, u_0(1)\big),
\mbox{ and } \|v \|_1 \le \eps  \Big\}.
\end{multline*}
 From now, we assume implicitly that $\| u_0\|_{C^1([0, 1])}$ is sufficiently small so that $\hD_\eps$ is not empty.
 For $v  \in \hD_\eps$, let $u = \cF(v)$ be the unique $C^1$-solution of the system 
\begin{equation}
\left\{\begin{array}{cl}
\partial_t u(t, x) = A(t, x) \partial_x u(t,x) + f \big(t, x, v(t, x) \big)& \mbox{ in } [0, T]  \times [0, 1],  \\[6pt]
u_{-}(t, 0)  = g \big(t, u_{+}(t, 0) \big)  & \mbox{ for } t \in [0, T],\\[6pt]
u_{+}(t, 1) = h \big(t,  v(t, \cdot) \big)& \mbox{ for } t \in [0,T], \\[6pt]
u(0, \cdot) = u_{0}(\cdot) &  \mbox{ in } [0, 1].
\end{array}\right.
\end{equation}
Here and in what follows, for notational ease, we ignore the dependence of $h$ on $u_0$ and denote $h(t, v(t, \cdot))$ instead of  $h(t, v(t, \cdot), u_0)$. As in the proof of \cite[Lemma 3.2]{CoronNg19} by  \eqref{lem-SL-cd-1} and \eqref{lem-SL-cd-2},  and the fact that $f$ and $g$ are of class $C^1$, one can prove that $\cF$ is contracting for $\|
\cdot \|_1$-norm provided that $L_2$ is large and $L_1$ is much larger than $L_2$. The condition $0 \le a<1$ and $1 \le p< + \infty$ are essential for the existence of $L_1$ and $L_2$.  \footnote{We here clarify a misleading point in
the definition of  $\cF(v)$ in \cite[(3.10)]{CoronNg19} in the proof of \cite[Lemma 3.2]{CoronNg19}. Concerning this definition,  in the RHS of \cite[(3.8)]{CoronNg19}, $v_{j+k}(t, 0)$ must be  understood as $(\cF (v))_{j+k}(t, 0)$ and $(\cF (v))_{j+k}(t, 0)$ is then determined by the RHS of \cite[(3.6) or (3.7)]{CoronNg19} as mentioned there. Related to this point,  $V_{j}(t, 0)$ for $ k+1 \le j \le k+m$ in \cite[(3.14)]{CoronNg19} and in the inequality just below must be replaced by $(\cF(v) - \cF(\hat v))_j$. The rest of the proof is unchanged.}
The existence and uniqueness of $u$ then follow. Moreover, there exist two constants $C_1, C_2>0$, independent of $u_0$ such that for  $\|  u_0\|_{C^1([0, 1])} \le C_1 \eps$ and $\|v\|_1 < \eps$, there exists a unique solution  $u \in \big(C^1([0, T] \times [0, 1])\big)^n$ and  moreover,
\begin{equation*}
\| u\|_{C^1([0, T] \times [0, 1])} \le C_2  \Big(\|  u_0\|_{C^1([0, 1])} + \sup_{t >0} \big( | h_1(t, u_0) | + | \partial_t h_1(t, u_0) | \big) \Big).
\end{equation*}
It follows from \eqref{lem-SL-cd-02} that for $\eps > 0$ small, there exists a constant $0< C_3 (\eps) < \eps$ small, independent of $u_0$,  such that for $\| u_0 \|_{C^1([0, 1])} \le C_3(\eps)$ and $v \in \hD_\eps$, then
\begin{equation}\label{lem-SL-bound}
\| \cF(v) \|_1 \le \eps  \mbox{ which implies in particular that } \cF(v) \in \hD_\eps.
\end{equation}
It is clear that $\cF(v) \in \hD_\eps$.

We claim that, for $\| u_0 \|_{C^1([0, 1])} \le C_3(\eps) $ and $\eps$ sufficiently small,
\begin{equation}\label{lem-SL-claim}
\mbox{$\cF$ is a contraction mapping w.r.t. $\| \cdot\|_1$ from $\hD_\eps$ into $\hD_\eps$.}
\end{equation}

Indeed, fix $\lambda \in (0, 1)$.  As in the proof of \cite[Lemma 3.2]{CoronNg19}, applying the characteristic method, and using \eqref{lem-SL-cd-1} and  \eqref{lem-SL-cd-2}, and the fact $f$ and $g$ are of class $C^1$, we obtain
\begin{equation}\label{lem-SL-p1}
\|\cF(\hat v) - \cF(v) \|_0 \le \lambda \|\hat v - v \|_1,
\end{equation}
if  $L_2$ is large and $L_1$ is much larger than $L_2$.  Set $U (t, x) = \partial_t u (t, x)$ for $(t, x) \in [0, T] \times [0, 1]$. We have
\begin{equation}
\left\{\begin{array}{cl}
\partial_t U(t, x) = A(t, x) \partial_x U(t,x)  +  \partial_t A(t, x) A(t, x)^{-1} U(t, x)
+ f_1 (t, x, v)   \mbox{ in } [0, T] \times [0, 1],  \\[6pt]
U_{-}(t, 0)  = g_1(t)   \mbox{ for } t \in [0, T],\\[6pt]
U_{+}(t, 1) = h_1 (t) \mbox{ for } t \in [0, T], \\[6pt]
U(0, x) = A(0, x) u_0'(x) + f\big(0, x, u_0(x) \big) \mbox{ in } [0, 1],
\end{array}\right.
\end{equation}
where
\begin{multline*}
f_1(t, x, v) = - \partial_t A(t, x) A^{-1}(t, x) f\big(t, x, v(t, x) \big) + \partial_t f \big(t, x, v(t, x) \big)  + \partial_{y} f\big(t, x, v(t, x) \big) \partial_t v(t, x).
\end{multline*}
$$
g_1(t) = \partial_t g(t, u_+(t, 0)) + \partial_{y_+} g(t, u_+(t, 0)) U_+(t, 0),
$$
$$
h_1(t) = \partial_t h(t, v(t, \cdot)) + \partial_y h (t, v(t, \cdot)) \partial_t v (t, \cdot).
$$
Note that, with $\hat u = \cF(\hat v)$ and $\hat U = \partial_t \hat u$,
\begin{multline*}
\Big|\partial_t g(t, \hat u_+(t, 0)) + \partial_{y_+} g(t, \hat u_+(t, 0)) \hat U_+(t, 0) \\[6pt]- \partial_t g(t, u_+(t, 0)) - \partial_{y_+} g(t, u_+(t, 0)) U_+(t, 0) \Big|  \\[6pt]
\mathop{\le}^{g \in C^2}
C \Big(  |\hat u_+(t, 0) - u_+(t, 0)| + |\hat U_+(t, 0) - U_+(t, 0)| \Big),
\end{multline*}
and
$$
|f_1(t, x,  \hat v) - f_1(t, x, v) | \mathop{\le}^{f \in C^2} C \Big( |\hat v(t, x) - v(t, x)| + |\partial_t \hat v(t, x) - \partial_t v(t, x)| \Big), 
$$
and by \eqref{lem-SL-cd-2} and \eqref{lem-SL-cd-3}, 
\begin{multline*}
\Big|\partial_t h(t, v(t, \cdot)) + \partial_y h (t, v(t, \cdot)) \partial_t v (t, \cdot) - \partial_t h(t, \hat v(t, \cdot)) - \partial_y h (t, \hat v(t, \cdot)) \partial_t \hat v (t, \cdot) \Big| \\[6pt]
\le C \Big( \| \big( (\hat v - v) (t, \cdot),  \partial_t (\hat v - v) (t, \cdot), \partial_x (\hat v - v) (t, \cdot) \big) \|_{C^0([0, a])} \\[6pt] +  \| \big( (\hat v - v) (t, \cdot),  \partial_t (\hat v - v) (t, \cdot), \partial_x (\hat v - v) (t, \cdot)  \big) \|_{L^p(0, 1)} \Big), 
\end{multline*}
if $\max\big\{ \|u  \|_{1}, \| v\|_1, \| \hat u \|_1, \|\hat v \|_1 \big\}  <  \eps_0$.  Again, as in the proof of \cite[Lemma 3.2]{CoronNg19}, applying  the characteristic method and using \eqref{lem-SL-cd-2} and \eqref{lem-SL-cd-3}, we also have, by \eqref{lem-SL-p1},
\begin{equation}\label{lem-SL-p2}
\|\partial_t \cF(\hat v) - \partial_t \cF(v) \|_0 \le \lambda \|\hat v - v \|_1.
\end{equation}
Since
$$
\partial_t (\hat u - u) (t, x) = A(t, x) \partial_x (\hat u - u)(t,x)  + f(t, x, \hat v(t, x)) - f(t, x, v(t, x)),
$$
and
$$
|f\big(t, x, \hat v(t, x) \big) - f\big(t, x, v(t, x) \big)| \le C |\hat v(t, x) - v(t, x)|,
$$
it follows from \eqref{lem-SL-p1} and \eqref{lem-SL-p2} that
\begin{equation*}
\|\cF(\hat v) -  \cF(v) \|_1 \le C \lambda \|\hat v - v \|_1.
\end{equation*}
Claim \eqref{lem-SL-claim} is proved.

The existence and uniqueness of solutions of \eqref{lem-SL-sys} in $\big(C^1([0, T] \times [0, 1])\big)^n$ now follow for $u_0$ satisfying $\| u_0 \|_{C^1([0, 1])} \le C_3(\eps)$. The proof is complete.
\end{proof}

We next establish the key result of this section. To this end, we first set, for $\tau>0$,
\begin{multline*}
\hat \D_\tau : =  \Big\{(\Xi, \varphi, w_0) \in   \big(C^1([0, + \infty) \times [0, 1]) \big)^n \times \big(C^1([0, 1])\big)^n \times \big(C^1([0, 1])\big)^n; \\[6pt]
\max \big\{\|\Xi\|_{C^1([0, + \infty) \times [0, 1])}, \| \varphi\|_{C^1([0, 1])}, \| w_0\|_{C^1([0, 1])} \big\}  <  \tau \Big\}
\end{multline*}
and, for $T>0$,
\begin{multline*}
\D_\tau : =  \Big\{(\Xi, u_0); (\Xi, 0, w_0) \in \hat \D_\tau; \Xi(0, \cdot) = w_0 (\cdot), \Xi(t, \cdot) = 0 \mbox{ for } t > T, \\[6pt]
\mbox{ and the compatibility conditions at $x = 0$ hold for the system \eqref{lem-QL-0} below}
 \Big\}.
\end{multline*}
The set $\D_\tau$ also depends on $T$ but we ignore this dependence explicitly for notational ease.

\medskip

We have

\begin{lemma} \label{lem-QL} Let  $T>0$, $f: [0, + \infty) \times [0, 1] \times \mR^n \to \mR^n$ be of class $C^2$ such that $f(t, x, 0) = 0$ for $(t, x ) \in [0, + \infty) \times [0, 1]$.
Assume that $B = \nabla \B (0) \in \cB$, $\Sigma$ is of class $C^2$,  and there exist $\tau > 0$ and
$$
H: [0, + \infty) \times \hat \D_\tau \to \mR^m
$$
such that $H$ is continuously differentiable w.r.t. $(t, \Xi, \varphi)$, and for some $C>0$,  $1 \le p <  + \infty$,  and $a \in [0, 1)$,  the following conditions hold, for $(\Xi, \varphi, w_0), (\hat \Xi, \hat  \varphi, w_0) \in \hat \D_\tau$ with $(\Xi, w_0), (\Xi, w_0)  \in \D_\tau$,
\begin{equation}\label{F-01}
H(t, \Xi,  \varphi, w_0) = H_1(t, w_0) + H_2(t, \Xi,  \varphi, w_0) \mbox{ with } H_1, \, H_2 \mbox{ are of class $C^1$}, 
\end{equation}
\begin{equation}\label{F-02}
\lim_{\eta \to 0} \sup_{\| u_0 \|_{C^1([0, 1])} \le \eta} \sup_{t >0} \Big( | H_1(t, u_0) | + | \partial_t H_1(t, u_0) | \Big) = 0, 
\end{equation}
\begin{equation}\label{F-1}
| H_2\big(t, \Xi,  \varphi, w_0 \big) | \le C \Big( \| (\varphi, \varphi') \|_{C^0([0, a])} +  \| (\varphi, \varphi') \|_{L^p(0, 1)} \Big),
\end{equation}
\begin{multline}\label{F-2}
| H \big(t, \hat \Xi ,  \hat \varphi, w_0  \big) - H \big(t, \Xi ,  \varphi, w_0 \big)| + | \partial_t H \big(t, \hat \Xi ,  \hat \varphi, w_0  \big) - \partial_t H \big(t, \Xi ,  \varphi, w_0 \big)|  \\[6pt]
 \le C \Big(\|\hat \Xi  - \Xi \|_{C^0([0, + \infty) \times [0, 1])} \| \hat \varphi\|_{C^1([0, 1])} + \| \hat \varphi - \varphi \|_{C^0([0, a])} +  \| \hat \varphi - \varphi \|_{L^p([0, 1])}\Big),
\end{multline}
\begin{equation}\label{Fv-1}
|\langle \partial_\varphi H(t, \Xi, \varphi, w_0), d\varphi \rangle|
\le C \Big( \| d\varphi \|_{C^0([0, a])} +  \| d\varphi \|_{L^p(0, 1)} \Big) \quad  \forall \,  d \varphi \in  \big(C^1([0, 1])\big)^n,
\end{equation}
\begin{equation}\label{Fxi-1}
\left|\frac{d}{dt} H(s, \Xi(t + \cdot, \cdot), \varphi, w_0) |_{s=t} \right|
\le C \Big(\|\Xi(t + \cdot, \cdot)\|_{C^1([0, + \infty) \times [0, 1])} + \| \varphi \|_{C^1([0, 1])} \Big)   \|\varphi \|_{C^1([0, 1])}.
\end{equation}
and, for $\eta > 0$ and for $0 \le |t' - t| \le \eta$, for $d \varphi, d \hat \varphi \in \big(C^1([0, 1])\big)^n$, 
\begin{multline}\label{Ftau-1}
\left| \frac{d}{ds} H_2 (s, \Xi(t' + \cdot, \cdot), \hat \varphi, w_0) |_{s=t'} -  \frac{d}{ds} H (s, \Xi(t + \cdot, \cdot), \varphi, w_0) |_{s=t}\right| \\[6pt]
+ \left|  \frac{d}{dt'} H_2 (s, \Xi(t' + \cdot, \cdot), \hat \varphi, w_0) |_{s=t'} - \frac{d}{dt} H (s, \Xi(t + \cdot, \cdot), \varphi, w_0) |_{s=t} \right| \\[6pt]
+ |  \langle \partial_\varphi H(t', \Xi, \hat \varphi, w_0), d\hat \varphi \rangle - \langle \partial_\varphi H(t, \Xi, \varphi, w_0), d\varphi \rangle| \\[6pt]
\le C \Big(\rho_1(c \eta, w_0) + \rho_2(c \eta, \varphi, \hat \varphi, d \varphi, d \hat \varphi)\Big),
\end{multline}
for some constant $c>0$ and some function $\rho_1$ such that 
$$
\lim_{\eta \to 0} \rho_1(\eta, w_0) =0, 
$$
where 
\begin{multline*}
\rho_2(\eta, \varphi, \hat \varphi, d \varphi, d \hat \varphi) = \| \mathop{\sup_{y}}_{|y - x| \le \eta}  \Big\{ |\varphi(y) - \hat \varphi(x)| +  |d\varphi(y) -  d \hat  \varphi(x)| \Big\} \|_{L^p(0, 1)} \\
+  \| \mathop{\sup_{y}}_{|y - x| \le \eta}  \Big\{  |\varphi(y) - \varphi(x)| +  |d\varphi(y) - d\varphi(x)| \Big\}\|_{C([0, a])}. 
\end{multline*}
Assume also  that for all $(\Xi, w_0)  \in \D_\tau$, the  system
\begin{equation}\label{lem-QL-0}
\left\{\begin{array}{cl}
\partial_t w(t, x) = \Sigma(x, \Xi (t, x)) \partial_x w(t,x)  + f(t, x, w(t, x))& \mbox{ in } [0, + \infty) \times [0, 1],  \\[6pt]
w_{-}(t, 0)  = \B \big( w_{+}(t, 0) \big) & \mbox{ for } t \in [0, + \infty),\\[6pt]
w_{+}(t, 1) = H \big(t, \Xi(t + \cdot, \cdot), w(t, \cdot), w_0 \big) & \mbox{ for } t \in [0, + \infty), \\[6pt]
w(0, \cdot) = w_0(\cdot) &  \mbox{ in } [0, 1]
\end{array}\right.
\end{equation}
has a unique $C^1$-solution satisfying $w(t, \cdot) = 0$ for $t > T$.  There exists $\eps > 0$ such that if $\|w (0, \cdot) \|_{C^1([0, 1])} < \eps$ and $w(0, \cdot)$ satisfies the compatibility conditions at $x=0$, then there is a unique solution $w \in \big(C^1([0, T] \times [0, 1])\big)^n$ of
\eqref{Sys-1} and  \eqref{bdry-w-0} with
\begin{equation}\label{wt1-conv}
w(t, 1) = H\big(t, w(t + \cdot, \cdot), w(t, \cdot), w_0 \big) \mbox{ for } t \in [0, + \infty).
\end{equation}
Moreover,
\begin{equation}\label{lem-QL-est}
\| w\|_{C^1([0, + \infty) \times [0, 1])} \le C \Big( \| w_0\|_{C^1([0, 1])} + \sup_{\| u_0 \|_{C^1([0, 1])} \le \eta} \sup_{t >0} \Big( | H_1(t, u_0) | + | \partial_t H_1(t, u_0) | \Big) \Big),
\end{equation}
for some positive constant independent of $w_0$ and $\eps$.
\end{lemma}

In \Cref{lem-QL} and what follows, $\Xi(t + \cdot, \cdot)$ denotes the function $(s, x) \mapsto \Xi(t + s,x)$ and $w(t + \cdot, \cdot)$ denotes the function $(s, x) \mapsto w(t + s, x)$.

The compatibility conditions at $x=0$ considered in the context of  \Cref{lem-QL} are
\begin{equation*}
w_{0, -}(0) = B\big(w_{0, +}(0) \big), 
\end{equation*}
and
\begin{multline*}\Big(\Sigma(0, \Xi (0, 0)) \partial_x w(0,0)  + f(0, 0, w(0, 0)) \Big)_- \\[6pt]
 = \nabla B (w_+(0, 0)) \Big( \Sigma(0, \Xi (0, 0)) \partial_x w(0,0)  + f(0, 0, w(0, 0)) \Big)_+.
\end{multline*}

The compatibility at $x=1$ of \eqref{lem-QL-0} is a part of the assumption of \Cref{lem-QL}.

Before giving the proof of \Cref{lem-QL}, let us discuss the motivation for the assumptions made. To this end,  we present one of its applications used in the proof of \Cref{thm1}. Consider the setting given in \Cref{thm1}; $f = 0$ then.  For $\Xi \in \big(C^1([0, + \infty) \times [0, 1])\big)^n$, define the flows
$$
\frac{d}{dt}x^\Xi_j(t, s, \xi) =   \lambda_j \Big(x^\Xi_j(t, s, \xi), \Xi \big(t, x^\Xi_j(t, s, \xi) \big) \Big) \quad \mbox{ and } \quad x_j^\Xi(s, s, \xi) = \xi  \mbox{ for } 1 \le j \le k,
$$
and
$$
\frac{d}{dt}x^\Xi_j(t, s, \xi) = -  \lambda_j \Big(x^\Xi_j(t, s, \xi), \Xi \big(t, x^\Xi_j(t, s, \xi) \big) \Big) \quad \mbox{ and } \quad x^\Xi_j(s, s, \xi) = \xi  \mbox{ for } k+1 \le j \le k + m.
$$
Here and in what follows, we only consider the flows with $x^\Xi_j(t, s, \xi) \in [0, 1]$ so that  $\Xi$ is well-defined.
Assume that $m>k$. Since $\nabla \B(0) \in \cB$, by the implicit theorem and the Gaussian elimination method,  there exist $M_k: U_k \to \mR$, \dots, $M_1: U_1 \to \mR$  of class $C^2$ for some neighborhoods $U_k$ of $0 \in \mR^{m-1}$, \dots, $U_1$ of $0 \in \mR^{m-k}$ such that, for $y_+ = (y_{k+1}, \cdots, y_{k +m})\tr \in \mR^m$ with sufficiently small norm,  the following facts hold
$$
\Big( \B(y_+) \Big)_k = 0 \mbox{ if } y_{k+m}  = M_k (y_{k+1}, \dots, y_{k+ m-1}),
$$
$$
\Big( \B(y_+) \Big) _k =  \Big(  \B(y_+) \Big)_{k-1} =  0 \mbox{ if } y_{k+m}  = M_k (y_{k+1}, \dots, y_{k+ m-1}),  y_{k+m-1}  = M_{k-1} (y_{k+1}, \dots, y_{k+ m-2}),
$$
\dots,
$$
\B(y_+) =  0 \mbox{ if } y_{k+m}  = M_k (y_{k+1}, \dots, y_{m+1}), \dots, y_{m+1}  = M_{1} (y_{k+1}, \dots, y_{m}).
$$
For $T> T_{opt}$, set $\delta = T - T_{opt}$.  Consider $\zeta_j$ and $\eta_j$ of class $C^1$ for $k+1 \le j \le k+m$ and for $t \ge 0$ satisfying
\begin{equation}\label{def-xi-eta-1}
\zeta_{j}(0) = w_{0, j}(1), \quad \zeta_{j}(t) = 0 \mbox{ for } t \ge \delta/2, \quad \eta_{j}(0) = 1, \quad \eta_{j}(t) = 0 \mbox{ for } t \ge \delta/2,
\end{equation}
and
\begin{equation}\label{def-xi-eta-2}
\zeta_j'(0) = \lambda_j\big(1, w_0(1)\big) w_{0, j}'(1), \quad \eta_j'(0)=0.
\end{equation}
For $(\Xi, \varphi, w_0) \in \D_\tau$ with small $\tau$,  set
\begin{multline}\label{bdry-1-claim}
\Big(H(t, \Xi, \varphi, w_0) \Big)_{m} =  \zeta_{k+m}(t)  \\[6pt]+ ( 1 - \eta_{k+m}(t)) M_{k} \Big(\varphi_{k+ 1}\big( x^\Xi_{k+1} (t, t+ t^\Xi_{m+k},  0) \big), \dots, \varphi_{k+ m-1}\big( x^\Xi_{k+ m - 1} (t, t+  t^\Xi_{m + k }, 0) \big)\Big),
\end{multline}
\begin{multline}\label{bdry-2-claim}
\Big(H(t, \Xi, \varphi, w_0)\Big)_{m-1} = \zeta_{k+m-1}(t)  \\[6pt]
+  ( 1-  \eta_{k+m-1}(t))  M_{k-1} \Big(\varphi_{k+ 1}\big(x^\Xi_{k+1} (t, t+ t^\Xi_{m + k - 1}, 0) \big), \dots, \varphi_{k+ m-2}\big(x^\Xi_{k+ m - 2} (t, t+  t^\Xi_{m + k -1}, 0) \big) \Big),
\end{multline}
\dots
\begin{multline}\label{bdry-m-claim}
\Big(H(t, \Xi, \varphi, w_0)\Big)_{m+1 - k} =  \zeta_{m+1}(t) \\[6pt]
+  (1- \eta_{m+1}(t))  M_{1} \Big(\varphi_{k+ 1} \big(x^\Xi_{k+1} (t, t + t^\Xi_{m + 1}, 0) \big), \dots, \varphi_{m}\big(x^\Xi_{m} (t, t+  t^\Xi_{m +1}, 0) \big)\Big),
\end{multline}
and
\begin{equation}\label{bdry-m-claim-1}
\Big(H(t, \Xi, \varphi, w_0)\Big)_{j} =  \zeta_{k+j}(t) \mbox{ for } 1 \le j \le m-k,
\end{equation}
where $t^\Xi_{j} = t^\Xi_{j}(t)$ are defined by
$$
x^\Xi_{m+k}(t+ t^\Xi_{m+k}, t, 1) = 0, \dots,  x^\Xi_{1+k}(t+ t^\Xi_{1+k}, t, 1) = 0 \mbox{ for } k+1 \le j \le k+ m.
$$
We now show that $H$ satisfies the assumptions given in \Cref{lem-QL} if $\|w_0\|_{C^1([0, 1])} \le \eps$ and $\eps$ is sufficiently small ($\tau$ is sufficiently small as well). We first note that the solutions of the system \eqref{lem-QL-0} are 0 for $t > T$ if $\| \Xi\|_{C^1([0, + \infty) \times [0, 1])}$ is sufficiently small. The proof of this fact follows from the choice of $M_j$ (see the proof of \eqref{nul-1-1}-\eqref{nul-m-1} in the proof of \Cref{thm1}).
One can easily check that \eqref{F-01}, \eqref{F-1}, \eqref{Fv-1}, and \eqref{Fxi-1} hold. Assertion \eqref{F-02} will be a consequence of our construction $\eta_j$ and $\zeta_j$ given later.  We are next concerned about \eqref{F-2}. It suffices to prove that
\begin{multline}\label{F-2-aux}
|H(t, \Xi, \varphi, w_0) - H(t, \hat \Xi, \varphi, w_0)| + |\partial_t H(t, \Xi, \varphi, w_0) - \partial_t H(t, \hat \Xi, \varphi, w_0)|  \\[6pt]
\le C \|\hat \Xi - \Xi  \|_{C^0([0, + \infty) \times [0, 1])} \| \varphi\|_{C^1([0, 1])}.
\end{multline}
We claim that, for $1 \le j \le k+m$.
\begin{equation}\label{claim-2-flow}
|x^{\hat \Xi}_{j}(t, s, \xi) - x^{\Xi}_j(t, s, \xi) | \le C \| \hat \Xi - \Xi\|_{C^0([0, + \infty) \times [0, 1])}
\end{equation}
for $(t, s, \xi)$ so that both flows are well-defined.
We only consider the case $k+1 \le j \le k+m$, the other cases can be proved similarly. We have
\begin{equation*}
|x^{\hat \Xi}_{j}(t, s, \xi) - x^{\Xi}_j(t, s, \xi) | \le C \| \hat \Xi - \Xi\|_{C^0([0, + \infty) \times [0, 1])} + C \int_{\min\{t, s \}}^{\max\{t, s \}}  |x^{\hat \Xi}_{j}(s', s, \xi) - x^{\Xi}_j(s', s, \xi) | \, ds'
\end{equation*}
and \eqref{claim-2-flow} follows.

Since, for $k+1 \le j \le k+m$,
$$
\int_{t}^{t  + t^{\hat \Xi}_{j} }  \lambda_j \Big(x^{\hat \Xi}_j(s, t, 1), \hat \Xi \big(t, x^{\hat \Xi}_j(s, t, 1) \big) \Big) \, ds  = 1 = \int_{t}^{t  + t^{\Xi}_{j} }  \lambda_j \Big(x^\Xi_j(s, t, 1), \Xi \big(t, x^\Xi_j(s, t, 1) \big) \Big) \, ds,
$$
it follows from \eqref{relation-lambda} and  \eqref{claim-2-flow} that
\begin{align}\label{3-flow}
|t^{\hat \Xi}_j - t^{\Xi}_j| \le & C \int_t^{t + \min\{t^{\hat \Xi}_j, t^{\Xi}_j \}} \Big( |x^{\hat \Xi}_j(s, t, 1) - x^\Xi_j(s, t, 1)| + \| \hat \Xi - \Xi\|_{C^0([0, + \infty) \times [0, 1])} \Big)  \, ds \\[6pt]
 \le & C  \| \hat \Xi - \Xi\|_{ C^0 ([0, + \infty) \times [0, 1]) }.  \nonumber
\end{align}
Combining \eqref{claim-2-flow} and \eqref{3-flow} yields \eqref{F-2-aux}. One can also verify \eqref{Ftau-1} by direct/similar  computations
and by using the fact 
$$
|x^{\Xi}_j(t', s', \xi')  - x^{\Xi}_j(t, s, \xi)| \le C \big(|t' - t| + |s' - s| + |\xi' - \xi|\big). 
$$

\medskip
We now give the

\begin{proof}[Proof of \Cref{lem-QL}] In what follows, for notational ease, we ignore the dependence of $H$ on $w_0$ and denote $H(t, \Xi,  \varphi(t, \cdot))$ instead of  $H(t, \Xi,  \varphi(t, \cdot), w_0)$.  Fix an appropriate $w^{(0)}$ such that $(w^{(0)}, w_0) \in \D_\tau$ and $\| w^{(0)}\|_{C^1([0, + \infty) \times [0, 1])} \le C \| w_0\|_{C^1([0, 1])}$; we thus assumed implicitly here that $\|w_0\|_{C^1([0, 1])}$ is sufficiently small.  For $l \ge 0$, let $w^{(l+1)}$  be the unique $C^1$-solution of
\begin{equation}
\left\{\begin{array}{cl}
\partial_t w^{(l+1)}(t, x) = \Sigma(x, w^{(l)} (t, x)) \partial_x w^{(l+1)}(t,x) +  f(t, x, w^{(l+1)}(t, x))& \mbox{ in } [0, + \infty) \times [0, 1],  \\[6pt]
w^{(l+1)}_{-}(t,  0)  = \B \big(w^{(l+1)}_{+}(t, 0) \big) & \mbox{ for } t \in [0, + \infty),\\[6pt]
w^{(l+1)}_{+}(t, 1) = H \big(t,  w^{(l)} (t + \cdot , \cdot), w^{(l+1)}(t, \cdot) \big) & \mbox{ for } t \in [0, + \infty), \\[6pt]
w^{(l+1)}(0, \cdot) = w_{0}(\cdot) &  \mbox{ in } [0, 1],
\end{array}\right.
\end{equation}
and set
$$
W^{(l)}(t, x) = \partial_t w^{(l)}(t, x)  \mbox{ for } (t, x)  \in [0,  + \infty) \times [0, 1].
$$
The existence and uniqueness of $w^{(l+1)}$ follows from \Cref{lem-SL}. Indeed, the compatibility conditions at $x=0$ follow from the fact $w^{(l)}(0,
\cdot) = w_0(\cdot)$ and the compatibility conditions at $x=1$ follow from the assumption on $H$ for the existence of $C^1$-solutions of the system  \eqref{lem-QL-0}.
We have
\begin{equation}\label{lem-QL-sysW}
\left\{\begin{array}{l}
\partial_t W^{(l+1)}(t, x) = \Sigma(x, w^{(l)} (t, x)) \partial_x W^{(l+1)} (t,x) \\[6pt] \qquad   \qquad \qquad \qquad + f_1(t, x) W^{(l+1)}(t, x)
+ f_2(t, x) \mbox{ for } (t, x) \in [0, + \infty) \times [0, 1],  \\[6pt]
W^{(l+1)}_{-}(t,  0)  = \nabla \B  \big(w^{(l+1)}(t, 0) \big) W^{(l+1)}_{+}(t, 0) \mbox{ for } t \in [0, + \infty),\\[6pt]
W^{(l+1)}_{+}(t, 1) = \partial_t H \big(t,  w^{(l)} (t + \cdot , \cdot), w^{(l+1)}(t, \cdot) \big) +  \langle \partial_{\Xi} H \big(t, w^{(l)}(t + \cdot, \cdot), w^{(l+1)}(t,  \cdot) \big), W^{(l)}(t + \cdot, \cdot) \rangle \\[6pt]
\qquad \qquad \qquad \qquad  \qquad \qquad + \langle \partial_\varphi H \big(t, w^{(l)}(t + \cdot, \cdot), w^{(l+1)}(t, \cdot) \big), W^{(l+1)} (t, \cdot) \rangle   \mbox{ for } t \in [0, + \infty), \\[6pt]
W^{(l+1)}(0, \cdot) = \Sigma(\cdot, w_0(x)) w_0'(\cdot)  + f(0, x , w_0(x)) \mbox{ in } [0, 1],
\end{array}\right.
\end{equation}
where
$$
f_1(t, x) = \partial_y \Sigma(x, w^{(l)}(t, x)) W^{(l)}(t,x)\Sigma^{-1}(x, w^{(l)} (t, x))  + \partial_y f(t, x, w^{(l+1)}(t, x)), 
$$
and
$$
f_2(t, x) = \partial_t f(t, x, w^{(l+1)}(t, x)) - \partial_y \Sigma(x, w^{(l)}(t, x)) W^{(l)}(t,x)\Sigma^{-1}(x, w^{(l)} (t, x))   f(t, x, w^{(l+1)}(t, x)).
$$
We have, since  $H_2 \big(t,  w^{(l)} (t + \cdot), 0 \big) = 0$ by \eqref{F-1},
\begin{multline}\label{lem-QL-nonlinear0}
| \partial_t H \big(t,  w^{(l)} (t + \cdot, \cdot), w^{(l+1)}(t, \cdot) \big)|
\mathop{\le}^{\eqref{F-01}, \eqref{F-2}} C \Big( |\partial_t H_1(t)| +  \| w^{(l)} \|_{C^{0}([0, + \infty) \times [0, 1])} \|w^{(l+1)} (t, \cdot)) \|_{C^1([0, 1])} \\[6pt] + \| (w^{(l+1) }(t, \cdot)\|_{C^0([0, a])} + \| (w^{(l+1)}(t, \cdot) \|_{L^p([0, 1])} \Big),
\end{multline}
\begin{multline}\label{lem-QL-nonlinear}
 \Big|\langle \partial_\Xi H \big(t, w^{(l)}(t + \cdot, \cdot), w^{(l+1)}(t,  \cdot) \big), W^{(l)}(t + \cdot, \cdot) \rangle \Big| \\[6pt]
\mathop{\le}^{\eqref{Fxi-1}} C \Big(  \|w^{(l)} \|_{C^1([0, + \infty) \times [0, 1])} +  \|w^{(l+1)}(t, \cdot) \|_{C^1([0, 1])} \Big) \|w^{(l+1)}(t, \cdot) \|_{C^1([0, 1])},
\end{multline}
and
\begin{multline*}
\Big| \langle \partial_\varphi H \big(t, w^{(l)}(t + \cdot, \cdot), w^{(l+1)}(t, \cdot) \big), W^{(l+1)} (t, \cdot) \rangle\Big| \\[6pt]
\mathop{\le}^{\eqref{Fv-1}} C \Big( \| W^{(l+1)}(t, \cdot) \|_{C^0([0, a])} +  \| W^{(l+1)}(t, \cdot) \|_{L^p(0, 1)} \Big).
\end{multline*}
By introducing $\| \cdot \|_0$ and $\| \cdot \|_1$ as in \eqref{0-norm} and \eqref{1-norm},  and using the above three inequalities, one can prove that
\begin{equation}\label{lem-WP-est-D}
\| w^{(l+1)} \|_{C^1([0, + \infty) \times [0, 1])} \le C \Big( \sup_{t >0} \big(|H_1(t)|  + |\partial_t H_1(t)| \big)+  \| w_0\|_{C^1([0, 1])}\Big),
\end{equation}
if   $\| w^{(l)} \|_{C^1([0, + \infty) \times [0, 1])} \le \eps$ and $\eps$ is sufficiently small. The smallness of $\eps$ is also used to absorb the second term of the RHS of \eqref{lem-QL-nonlinear0} and  the RHS of \eqref{lem-QL-nonlinear}.
It follows  from \eqref{F-02} that there exists a constant $0< C_3(\eps)< \eps$, independent of $w_0$ such that
\begin{equation}\label{lem-QL-C1}
\| w^{(l)}\|_{C^1([0, + \infty) \times [0, 1])} \le C \eps,
\end{equation}
if
$$
\|w_0\|_{C^1[0, 1]} \le C_3(\eps) \mbox{ and } \eps \mbox{ is sufficiently small}.
$$
This fact will be assumed from now on.

Set, for $l \ge 1$,
$$
V^{(l)} = w^{(l)} - w^{(l-1)} \mbox{ in } [0, +\infty) \times [0, 1].
$$
We have
\begin{equation*}
\left\{\begin{array}{l}
\partial_t V^{(l+1)}(t, x) = \Sigma(x, w^{(l)} (t, x)) \partial_x V^{(l+1)}(t,x) \\[6pt]
\qquad \qquad \qquad \qquad
+ \Big( \Sigma(x, w^{(l)} (t, x)) - \Sigma(x, w^{(l-1)} (t, x))  \Big) \partial_x w^{(l)}(t, x)   \\[6pt]
\qquad \qquad \qquad \qquad \qquad + f(t, x, w^{(l+1)}(t, x)) - f(t, x, w^{(l)}(t, x))
 \mbox{ in } [0, + \infty) \times [0, 1],  \\[6pt]
V^{(l+1)}_{-}(t, 0)  = \B \big(w^{(l+1)}_{+}(t, 0) \big) - \B \big(w^{(l)}_{+}(t, 0) \big)  \mbox{ for } t \in [0, + \infty),\\[6pt]
V^{(l+1)}_{+}(t, 1) = H \big(t,  w^{(l)} (t + \cdot , \cdot), w^{(l+1)}(t, \cdot) \big)  - H \big(t,  w^{(l-1)} (t + \cdot , \cdot), w^{(l)}(t, \cdot) \big)  \mbox{ for } t \in [0, + \infty), \\[6pt]
V^{(l+1)}(0, \cdot) = 0  \mbox{ in } [0, 1].
\end{array}\right.
\end{equation*}
Note that, by  \eqref{lem-QL-C1},
\begin{equation*}
 \Big| \Big( \Sigma(x, w^{(l)} (t, x)) - \Sigma(x, w^{(l-1)} (t, x))  \Big) \partial_x w^{(l-1)}(t, x) \Big| \mathop{\le}^{\Sigma \in C^1} C \eps |V^{(l)} (t, x)|,
\end{equation*}
\begin{equation*}
| f(t, x, w^{(l+1)}(t, x)) - f(t, x, w^{(l)}(t, x))| \mathop{\le}^{f \in C^1} C |V^{(l+1)}(t, x)|, 
\end{equation*}
\begin{equation*}
| \B \big(w^{(l+1)}_{+}(t, 0) \big) - \B \big(w^{(l)}_{+}(t, 0) \big)| \mathop{\le}^{\B \in C^1} C |V^{(l+1)}_+(t, 0)|,
\end{equation*}
\begin{multline*}
\Big| H \big(t, w^{(l)} (t + \cdot , \cdot), w^{(l+1)}(t, \cdot) \big)
- H \big(t,  w^{(l-1)} (t + \cdot , \cdot), w^{(l)}(t, \cdot) \big) \Big| \\[6pt]
\mathop{\le}^{\eqref{F-2}} C  \Big(\eps  \| V^{(l)} (t +  \cdot , \cdot) \|_{C^0([0, +\infty] \times [0, 1])}  +  \| V^{(l+1)}(t, \cdot) \|_{C^0([0, a])} + \| V^{(l+1)}(t, \cdot) \|_{L^p([0, 1])} \Big).
\end{multline*}
Set
$$
Y_l(t) = \max_{1 \le i \le n} \mathop{\max}_{(s, x) \in [0, t] \times [0, 1]} | e^{-L_1 s - L_2 x} V^{(l)}_{ i} (s, x)|.
$$
It follows that, provided that $L_2$ is large and $L_1$ is much larger than $L_2$,
\begin{equation*}
Y_{l+1}(t) \le \int_0^t \Big( \alpha Y_{l+1} (s) + \beta Y_{l}(s)   \Big) \, ds + C \eps Y_l(T),
\end{equation*}
for some $\alpha, \beta > 0$. By multiplying the above inequality with $e^{-Lt}$ for some large positive constant $L$, one can derive that, for $\eps$ sufficiently small, 
$$
\max_{t \in [0, T]} Y_{l+1}(t) e^{-L t} \le \frac{1}{2} \max_{t \in [0, T]} Y_{l}(t) e^{-L t}.
$$
This implies
\begin{equation}\label{lem-QS-m1}
w^{(l)} \mbox{ converges in } C^0([0, + \infty) \times [0, 1]).
\end{equation}

Set
$$
\rho(\eta, w^{(l)}) = \sup_{t, x}  e^{- L_1 t -  L_2 x}\mathop{\mathop{\sup}_{t', x'}}_{|(t,x) - (t', x')| \le \eta} \Big| \Big( \partial_t \big(w^{(l)}(t', x') - w^{(l)}(t, x) \big) , \partial_x \big(w^{(l)}(t', x') - w^{(l)}(t, x) \big) \Big) \Big|
$$
and
$$
\rho(\eta, w_0) = \sup_{|x - x'| \le \eta} | w_0'(x') - w_0'(x)|.
$$
Define the flows
$$
\frac{d}{dt}x^{(l)}_j(t, s, \xi) =   \lambda_j \Big(x^{(l)}_j(t, s, \xi), w^{(l)} \big(t, x^{(l)}_j(t, s, \xi) \big) \Big) \quad \mbox{ and } \quad x_j^{(l)}(s, s, \xi) = \xi  \mbox{ for } 1 \le j \le k,
$$
and
$$
\frac{d}{dt}x^{(l)}_j(t, s, \xi) = -  \lambda_j \Big(x^{(l)}_j(t, s, \xi), w^{(l)} \big(t, x^{(l)}_j(t, s, \xi) \big) \Big) \quad \mbox{ and } \quad x^{(l)}_j(s, s, \xi) = \xi  \mbox{ for } k+1 \le j \le k + m.
$$
By \eqref{relation-lambda} and the fact $\|w^{(l)} \|_{C^1([0, + \infty) \times [0, 1])} \le C \eps$, one has 
\begin{equation}\label{lem-QS-flow}
|x^{(l)}_j(t', s', \xi')  - x^{(l)}_j(t, s, \xi)| \le C \big(|t' - t| + |s - s'| + |\xi' - \xi|\big). 
\end{equation} 
Using \eqref{Ftau-1} and \eqref{lem-QS-flow}, and considering \eqref{lem-QL-sysW}, one can prove that
\begin{equation}\label{lem-QS-m2}
\rho(\eta, w^{(l)}) \le C \rho(C \eta, w_0) + C \eta + C \rho_1(C \eta, w_0).
\end{equation}

Combining \eqref{lem-QS-m1} and \eqref{lem-QS-m2}, and applying the Ascoli theorem, one derives that
$$
w^{(l)} \mbox{ converges in } \big(C^1([0, + \infty) \times [0, 1])\big)^n.
$$
It is clear that the limit is a $C^1$-solution of
\eqref{Sys-1}, \eqref{bdry-w-0}, and \eqref{wt1-conv}.

We next establish the uniqueness. Assume that $w$ and $\hat w$ are two $C^1$-solutions of \eqref{Sys-1}, \eqref{bdry-w-0}, and \eqref{wt1-conv}. Set $u = \hat w - w$ in $[0, + \infty) \times [0, 1]$.
Then
\begin{equation*}
\partial_t u (t, x) = A(t, x) \partial_x u (t, x) + \tilde f(t, x, u(t, x)),
\end{equation*}
where
$$
A(t, x) = \Sigma(x, w(t, x)),
$$
\begin{multline*}
\tilde f(t, x, u(t, x)) = \Big( \Sigma\big(x, w(t, x) + u(t, x)\big) - \Sigma\big(x, w(t, x) \big) \Big) \partial_x \hat w(t, x) \\[6pt]
+ f\big(t, x, w(t, x) + u(t, x) \big) - f\big(t, x, w(t, x)\big).
\end{multline*}
Moreover,
$$
u_-(t, 0) = g(t, u_+(t, 0)): = \B(w_+(t,0) + u_+(t, 0)) - \B(w_+(t, 0)),
$$
$$
u_+(t, 0) = h(t, u(t + \cdot, \cdot)):= H(t, w+ u, w+ u) - H(t, w, w),
$$
and
$$
u(t=0, \cdot) = 0.
$$
Note that
$$
|\tilde f(t, x, u(t, x))| \le C |u(t, x)|,
$$
$$
|g(t, u_+(t,0))| \le C |u_+(t, 0)|,
$$
and
$$
|h(t, u(t + \cdot, \cdot))| \mathop{\le}^{\eqref{F-2}} C \Big( \eps \| u \|_{C^0([0, + \infty) \times [0, 1])} +  \| u(t, \cdot) \|_{C^0([0, a])} +  \| u(t, \cdot) \|_{L^p([0, 1])} \Big).
$$
Let $U \in \big(C([0, + \infty) \times [0, 1])\big)^n$, with $U(t, \cdot) = 0$ for $t > T$,  be a  solution of the system
\begin{equation*}
\left\{\begin{array}{cl}
\partial_t U (t, x) - A(t, x) \partial_x U (t, x) = \tilde f(t, x, u(t, x)) & \mbox{ in } [0, + \infty) \times [0, 1], \\[6pt]
U_-(t, 0) = g \big(t, U_+(t, 0) \big) & \mbox{ for } t \in [0, + \infty),  \\[6pt]
U_+(t, 0) = h \big(t, u(t + \cdot, \cdot) \big) &  \mbox{ for } t \in [0, + \infty),  \\[6pt]
U(t=0, \cdot) = 0 & \mbox{ in } [0, 1],
\end{array}\right.
\end{equation*}
and set
$$
Y(t) = \max_{1 \le i \le n} \mathop{\max}_{(s, x) \in [0, t] \times [0, 1]} | e^{-L_1 s - L_2 x} U_i(s, x)|
$$
and
$$
Z(t) = \max_{1 \le i \le n} \mathop{\max}_{(s, x) \in [0, t] \times [0, 1]} | e^{-L_1 s - L_2 x} u_i(s, x)|.
$$
As in the proof of \cite[Lemma 3.2]{CoronNg19}, one can prove that, if  $L_2$ is large and $L_1$ is much larger than $L_2$,
$$
Y(t) \le C \int_0^t (Y(s) + Z(s)) \, ds  + C \eps Z(T).
$$
By multiplying 
the above inequality with $e^{-Lt}$, for some large positive constant $L$, one has 
$$
\max_{t \in [0, T]} Y(t) e^{-L t} \le \frac{1}{2} \max_{t \in [0, T]} Z(t) e^{-L t}.
$$
if $\eps$ is sufficiently small.
As a consequence, by taking $U = u$,
one has, for $\eps$ sufficiently small, $$
u=0
$$
and the uniqueness follows.  The proof is complete.
\end{proof}

\begin{remark} \rm The proof of \Cref{lem-QL} is inspired from \cite{CoronNg19} using the approach for quasilinear hyperbolic equations in \cite[Chapter 1]{LY85} and \cite[Chapter 3]{Bressan00}.

\end{remark}

\subsection{Proof of \Cref{thm1}} We consider two cases $m>k$ and $m \le k $ separately.

\medskip
\noindent {\bf Case 1: $m> k$.}  Consider the last equation of \eqref{bdry-w-0}. Impose the condition $w_{k}(t, 0) = 0$. Using \eqref{cond-B-1} with $i = 1$ and the implicit function theorem, one can then write the last  equation of \eqref{bdry-w-0} under the form
\begin{equation}\label{kd-1}
w_{m+k}(t, 0) = M_k  \Big( w_{k +1 }(t, 0), \cdots, w_{m + k - 1} (t, 0) \Big),
\end{equation}
for some $C^2$ nonlinear map  $M_k$ from $U_k$ into $\mR$ for some neighborhood $U_k$ of $0 \in \mR^{m-1}$ with $M_k(0) = 0$ provided that $|w_{+}(t, 0)|$ is sufficiently small.

Consider the last two equations of \eqref{bdry-w-0} and impose the condition $w_{k}(t, 0) = w_{k-1}(t, 0) = 0$. Using \eqref{cond-B-1} with $i = 2$ and  the Gaussian elimination approach,  one can then write these two equations under the form \eqref{kd-1} and
\begin{equation}\label{kd-2}
w_{m+k-1}(t, 0) = M_{k-1} \Big( w_{k +1 }(t, 0), \cdots, w_{m + k - 2} (t, 0) \Big),
\end{equation}
for some $C^2$ nonlinear map  $M_{k-1}$ from $U_{k-1}$ into $\mR$ for some neighborhood $U_{k-1}$ of $0 \in \mR^{m-2}$ with $M_{k-1}(0) = 0$ provided that $|w_{+}(t, 0)|$ is sufficiently small, etc. Finally, consider the  $k$ equations of \eqref{bdry-w-0} and impose the condition $w_{k}(t, 0) = \dots = w_{1}(t, 0) = 0$. Using \eqref{cond-B-1} with $i = k$ and the Gaussian elimination approach,  one can then write these $k$ equations under the form \eqref{kd-1}, \eqref{kd-2}, \dots, and \begin{equation}\label{kd-k}
w_{m+1}(t, 0) = M_{1} \Big(w_{k +1 }(t, 0), \cdots, w_{m} (t, 0) \Big),
\end{equation}
for some $C^2$ nonlinear map  $M_1$ from $U_1$ into $\mR$ for some neighborhood $U_1$ of $0 \in \mR^{m-k}$ with $M_1(0) = 0$ provided that $|w_{+}(t, 0)|$ is sufficiently small. These nonlinear maps  $M_1, \dots, M_k$ will be used in the construction of feedbacks.

We next introduce the flows along the characteristic curves. Set
$$
\frac{d}{dt}x_j(t, s, \xi) =   \lambda_j \Big(x_j(t, s, \xi), w \big(t, x_j(t, s, \xi) \big) \Big) \quad \mbox{ and } \quad x_j(s, s, \xi) = \xi  \mbox{ for } 1 \le j \le k,
$$
and
$$
\frac{d}{dt}x_j(t, s, \xi) = -  \lambda_j \Big(x_j(t, s, \xi), w \big(t, x_j(t, s, \xi) \big) \Big) \quad \mbox{ and } \quad x_j(s, s, \xi) = \xi  \mbox{ for } k+1 \le j \le k + m.
$$
We do not precise at this stage the domain of the definition of $x_j$. Later, we only consider the flows in the regions where the solution $w$ is well-defined.

To arrange the compatibility of our controls, we introduce auxiliary variables
satisfying autonomous dynamics, which  will be defined later. Set $\delta = T- T_{opt} > 0$. For $t \ge 0$, define, for $k+1 \le j \le k+m$,
\begin{equation}
\zeta_{j}(0) = w_{0, j}(1), \quad \zeta_j'(0) = \lambda_j\big(0, w_0(1)\big) w_{0, j}'(1), \quad \zeta_{j}(t) = 0 \mbox{ for } t \ge \delta/2,
\end{equation}
and
\begin{equation}
\eta_j(0) = 1, \quad \eta_j'(0)=0, \quad \eta_{j}(t) = 0 \mbox{ for } t \ge \delta/2.
\end{equation}
We will construct the dynamics for $\zeta_j$ and $\eta_j$  at the end of the proof of \Cref{thm1}.

\medskip
We are ready to construct a feedback law leading to finite-time stabilization in the time $T$.  Let $t_{m+k}$ be such that
$$
x_{m+k}(t+ t_{m+k}, t, 1) = 0.
$$
It is clear that $t_{m+k}$ depends only on the current state $w(t, \cdot)$. Let $D_{m+k} = D_{m+k}(t) \subset \mR^2$ be the open set whose boundary is $\{ t \} \times [0, 1]$, $[t, t+ t_{m+k}] \times \{0\}$, and $\Big\{ (s, x_{m+k}(s,  t, 1)); \ s \in [t, t+ t_{m+k}] \Big\}$. Then $D_{m+k}$ depends only on the current state as well. This implies
$$
x_{k+1}(t, t+ t_{m+k}, 0), \dots, x_{k+m-1}(t, t+ t_{m+k}, 0) \mbox{ are well-defined by the current state $w(t, \cdot)$}.
$$
As a consequence, the feedback
\begin{multline}\label{bdry-1}
w_{m+ k}(t, 1) =  \zeta_{m+k}(t)  \\[6pt]+ ( 1 - \eta_{m+k}(t)) M_{k} \Big(w_{k+ 1}\big(t, x_{k+1} (t, t+ t_{m+k},  0) \big), \dots, w_{k+ m-1}\big(t, x_{k+ m - 1} (t, t+  t_{m + k }, 0) \big)\Big)
\end{multline}
is well-defined by the current state $w(t, \cdot)$.

We then consider the system \eqref{Sys-1}, \eqref{bdry-w-0}, and the feedback \eqref{bdry-1}. Let $t_{m+k-1}$ be such that
$$
x_{m+k-1}(t+ t_{m+k-1}, t, 1) = 0.
$$
It is clear that $t_{m+k-1}$ depends only on the current state $w(t, \cdot)$ and the feedback law \eqref{bdry-1}. Let $D_{m+k-1} = D_{m+k-1}(t) \subset \mR^2$ be the open set whose boundary is $\{ t \} \times [0, 1]$, $[t, t+ t_{m+k-1}] \times \{0\}$, and $\Big\{ (s, x_{m+k-1}(s, t, 1)); \ s \in [t, t+  t_{m+k-1}] \Big\}$. Then $D_{m+k-1}$  depends only on the current state. This implies
$$
x_{k+1}(t, t+ t_{m+k-1}, 0), \dots, x_{k+m-2}(t, t+  t_{m+k-1}, 0) \mbox{ are well-defined by the current state $w(t, \cdot)$}.
$$
As a consequence, the feedback
\begin{multline}\label{bdry-2}
w_{m+ k - 1}(t, 1) = \zeta_{m+k-1}(t)  \\[6pt]
+  ( 1- \eta_{m+k-1}(t)) M_{k-1} \Big(w_{k+ 1}\big(t, x_{k+1} (t, t+ t_{m + k - 1}, 0) \big), \dots, w_{k+ m-2}\big(t, x_{k+ m - 2} (t, t+  t_{m + k -1}, 0) \big) \Big)
\end{multline}
is well-defined by the current state $w(t, \cdot)$.

We continue this process and finally reach the system \eqref{Sys-1}, \eqref{bdry-w-0}, \eqref{bdry-1}, \dots
\begin{multline}\label{bdry-m-1}
w_{m+ 2}(t, 1) = \zeta_{m+2}(t)  \\[6pt]
+  (1 -  \eta_{m+2}(t)) M_{2} \Big(w_{k+ 1}\big(t, x_{k+1} (t, t+  t_{m + 2}, 0) \big), \dots, w_{m+1}\big(t, x_{m+1} (t, t+  t_{m + 2}, 0) \big) \Big).
\end{multline}
Let $t_{m+1}$ be such that
$$
x_{m+1}(t+ t_{m+1}, t, 1) = 0.
$$
It is clear that $t_{m+1}$ depends only on the current state $w(t, \cdot)$ and the feedback law \eqref{bdry-1}, \dots, \eqref{bdry-m-1}. Let $D_{m+1} = D_{m+1}(t) \subset \mR^2$ be the open set whose boundary is $\{ t \} \times [0, 1]$, $[t, t+ t_{m+1}] \times \{0\}$, and $\Big\{ (s, x_{m+1}(s, t, 1)); \ s \in [t, t + t_{m+1}] \Big\}$. Then $D_{m+1}$ depends only on the current state. This implies
$$
x_{k+1}(t, t+ t_{m+1}, 0), \dots, x_{m}(t, t+ t_{m+1}, 0) \mbox{ are well-defined by the current state $w(t, \cdot)$}.
$$
As a consequence, the feedback
\begin{multline}\label{bdry-m}
w_{m+ 1}(t, 1) =  \zeta_{m+1}(t) \\[6pt]
+ (1 -  \eta_{m+1}(t))  M_{1} \Big(w_{k+ 1} \big(t, x_{k+1} (t, t + t_{m + 1}, 0) \big), \dots, w_{m}\big(t, x_{m} (t, t+  t_{m +1}, 0) \big)\Big)
\end{multline}
is well-defined by the current state $w(t, \cdot)$.

To complete the feedback for the system, we consider, for $k+1 \le j \le m$,
\begin{equation}\label{bdry-k-m}
w_{j}(t, 1) = \zeta_j(t),
\end{equation}

We will establish that  the feedback constructed gives the  finite-time stabilization in the time $T$ if $\eps$ is sufficiently small. To this end, we first claim that
\begin{equation}\label{thm1-claim-1}
\mbox{the system \eqref{Sys-1}, \eqref{bdry-w-0}, \eqref{bdry-1}, \dots, \eqref{bdry-m} is well-posed if $\eps$ is sufficiently small}.
\end{equation}
Indeed, it is clear to see that the feedback is given by
$$
H(t, w(t+ \cdot), w(t, \cdot),w_0),
$$
where $H$ is given by \eqref{bdry-1-claim}-\eqref{bdry-m-claim-1}.  The well-posedness for the feedback law is now a consequence of \Cref{lem-QL}
through  the example mentioned and examined right after it.

From \eqref{def-xi-eta-1} and \eqref{def-xi-eta-2}, we have, for $t \ge \delta/2$,
$$
\zeta_j (t) = 0 \mbox{ for } k+1 \le j \le k+m.
$$
It follows that, for $t \ge \delta/2$, the feedback law \eqref{bdry-1}, \dots, \eqref{bdry-m} has the form
\begin{equation}\label{bdry-1-t0}
w_{m+ k}(t, 1) =  M_{k} \Big(w_{k+ 1}\big(t, x_{k+1} (t, t+  t_{m+k},  0) \big), \dots, w_{k+ m-1}\big(t, x_{k+ m - 1} (t, t+  t_{m + k }, 0) \big)\Big),
\end{equation}
\begin{equation}\label{bdry-2-t0}
w_{m+ k - 1}(t, 1) = M_{k-1} \Big(w_{k+ 1}\big(t, x_{k+1} (t, t+  t_{m + k - 1}, 0) \big), \dots, w_{k+ m-2}\big(t, x_{k+ m - 2} (t, t+  t_{m + k -1}, 0) \big) \Big),
\end{equation}
\dots
\begin{equation}\label{bdry-m-t0}
w_{m+ 1}(t, 1) =  M_{1} \Big(w_{k+ 1} \big(t, x_{k+1} (t, t+  t_{m + 1}, 0) \big), \dots, w_{m}\big(t, x_{m} (t, t+  t_{m +1}, 0) \big)\Big).
\end{equation}

Set
$$
\hat t = \max\{\hat t_{k+1}, \dots, \hat t_{k+ m} \},
$$
where $\hat t_j$, for $k+1 \le j \le k+ m $, is defined by
$$
x_{j}(\hat t_j + \delta/2, \delta/2, 1) = 0.
$$
It follows from the characteristic method that
\begin{multline*}
w_{j}(t, \cdot) = 0 \mbox{ for } t \ge \hat t + \delta/2  \mbox{ for } k+1 \le j \le m, \\[6pt]
\mbox{ then for } j= m+1, \mbox{ then for } j =m+2, \dots, \mbox{ then for } j = m+k.
\end{multline*}
Using the characteristic method again, we have, by the choice of $M_k$,
\begin{equation}\label{nul-1}
w_{k}(t, 0) = 0 \mbox{ for } t \ge \delta/2 +  \hat t_{m+k},
\end{equation}
 by the choice of $M_k$ and $M_{k-1}$,
\begin{equation}\label{nul-2}
w_{k-1}(t, 0) = 0 \mbox{ for } t \ge \delta/2 + \hat t_{m+k-1},
\end{equation}
\dots, and,  by the choice of $M_k$, $M_{k-1}$, \dots, $M_1$,
\begin{equation}\label{nul-3}
w_{1}(t, 0) = 0 \mbox{ for } t \ge \delta/2 + \hat t_{m+1}.
\end{equation}

Let $\hat t_{k}$, \dots, $\hat t_1$ be such that
\begin{equation}\label{time-1-1}
x_{k}(\hat t_k + \delta/2 + \hat t_{m+k},   \delta/2 + \hat t_{m+k}, 0) = 1,
\end{equation}
\dots,
\begin{equation}\label{time-m-1}
x_{1}(\hat t_1 + \delta/2 + \hat t_{m+1},   \delta/2 + \hat t_{m+1}, 0) = 1.
\end{equation}
Using the characteristic method, we derive that
\begin{equation}\label{nul-1-1}
w_k(t, \cdot) = 0 \mbox{ for } t \ge \delta/2 + \hat t_{m+k} + \hat t_{k},
\end{equation}
\dots,
\begin{equation}\label{nul-m-1}
w_1(t, \cdot) = 0 \mbox{ for } t \ge \delta/2 + \hat t_{m+1} + \hat t_{1}.
\end{equation}

The conclusion follows by noting that
$$
|\hat t_j - \tau_j| \le \delta/ 4 \mbox{ for } 1 \le j \le k+ m,
$$
if $\eps$ is sufficiently small thanks to \eqref{F-02}  and  \eqref{lem-QL-est}.

\medskip
\noindent {\bf Case 2: $m \le k$.}  We consider the following feedback law
\begin{multline*}
w_{m+ k}(t, 1) =  \zeta_{m+k}(t)  \\[6pt]+  (1 - \eta_{m+k}(t)) M_{k} \Big(w_{k+ 1}\big(t, x_{k+1} (t, t+ t_{m+k},  0) \big), \dots, w_{k+ m-1}\big(t, x_{k+ m - 1} (t, t+  t_{m + k }, 0) \big)\Big),
\end{multline*}
\dots
\begin{equation*}
w_{k+2}(t, 1) = \zeta_{k+2}(t)  \\[6pt]
+  ( 1-\eta_{k+2}(t))  M_{2} \Big(w_{k+ 1}\big(t, x_{k+1} (t, t+ t_{k+2}, 0) \big) \Big),
\end{equation*}
and
$$
w_{k+1}(t, 1) = \zeta_{k+1}(t).
$$
The conclusion now follows by the same arguments. The details are omitted.

\medskip It remains to construct a dynamics for $\zeta_j$ and $\eta_j$. To this end, inspired by \cite{CVKB13, 2014-Perrollaz-Rosier}, we write $\zeta_j = \varphi_j + \psi_j$ where $\varphi_j$ and $\psi_j$ satisfy the dynamics
\begin{equation}\label{eq-vpj}
\varphi_j'(t) = - \frac{\alpha \varphi_j}{(\varphi_j^2 + \psi_j^2)^{1/3}}
\quad \mbox{ and } \quad
\psi_j'(t) = - \frac{\beta \varphi_j}{(\varphi_j^2 + \psi_j^2)^{1/3}},
\end{equation}
with $Y = (\varphi_j(0)^2 + \psi_j(0)^2)^{1/3}$,
\begin{equation}\label{cd-vpj}
\varphi_j(0) + \psi_j(0) = a,  \quad - \alpha \varphi_j(0) - \beta \psi_j (0) = b Y,
\end{equation}
where $a = w_{0, j}(0)$ and $b= \lambda_j (0, w_0(1)) w_{0, j}'(1)$. Here $\alpha$ and $\beta$ are two {\it distinct} real numbers.
We now show that under appropriate choice of $\alpha$ and $\beta$, $\varphi_j(0)$ and $\psi_j(0)$ can be chosen as  continuous functions of $a$ and $b$ for $|(a, b)|$ sufficiently small. Indeed, consider the
equation $P_{a, b}(Y) = 0$, where
\begin{equation}
P_{a, b}(Y): = (\alpha - \beta)^2 Y^3 - \Big( 2 b^2 Y^2 + 2 ab(\alpha + \beta)Y + a^2 (\alpha^2 + \beta^2) \Big).
\end{equation}
One has,  for $Y>0$ and $P_{a, b}(Y) =0$,
$$
Y P_{a, b}'(Y) = 2 b^2 Y^2 + 4 ab (\alpha + \beta) Y + 3 (\alpha^2 + \beta^2)a^2.
$$
In particular,
$$
P_{a, b}'(Y) > 0 \mbox{ if } \alpha^2 + \beta^2 - 4 \alpha \beta >0 \mbox{ and if }  ab \neq 0, 
$$
and the equation $P_{a, b}(Y) = 0$ has a unique positive solution in this case.
In the case $ab=0$ and $a^2 + b^2 >0$, there is a unique positive solution of $P_{a, b}(Y) = 0$ and in the case $a=b =0$, there is a unique solution $Y=0$. Fix $\alpha$ and $\beta$ such that $\alpha^2 + \beta^2 - 4 \alpha \beta \neq 0$ and $\alpha \neq \beta$. Denote $\bar Y(a, b)$  the unique positive solution in the case $a^2 + b^2 >0$ and 0 for $(a, b) = (0, 0)$. It suffices to prove that $\bar Y(a, b)$ is continuous with respect to $(a, b)$ for small $|(a, b)|$.
Since $P_{a, b}(1)>0$ if $|(a, b)|$ is sufficiently small and $P_{a,b}(0) <  0$ if $a \neq 0$, it follows that $\bar Y$ is bounded in a neighborhood $O$ of $(0, 0)$. Since $P_{a, b}(Y) = 0$ has a unique non-negative solution for $a \neq 0$, it follows that $\bar Y$ is continuous in $O \setminus \{ (a, b);  a =0\}$. Since $ \alpha^2 + \beta^2 - 4 \alpha \beta >0$, one has
$$
\frac{3}{2} b^2 Y^2 + 2 ab(\alpha + \beta)Y + a^2 (\alpha^2 + \beta^2)  \ge 0.
$$
It follows that
$$
P_{a, b}(Y) \le (\alpha - \beta)^2 Y^3 - \frac{1}{2} b^2 Y^2.
$$
This implies the continuity of $\bar Y$ on
$O \cap \{ (a, b);  a =0 \mbox{ and } b \neq 0\}$.
 The continuity of $\bar Y$ at $(0, 0)$ is a consequence of the fact $P_{0, 0}(Y) = 0$ implies $Y = 0$.

Similarly, one can build the dynamics for $\eta_j$. We now have $a=1$ and $b=0$. we write $\eta_j = \widetilde \varphi_j + \widetilde  \psi_j$ where $\widetilde  \varphi_j$ and $\widetilde \psi_j$ satisfy the dynamics
$$
\widetilde  \varphi_j'(t) = - \frac{\lambda^{5/3} \alpha \widetilde \varphi_j}{(\widetilde  \varphi_j^2 + \widetilde  \psi_j^2)^{1/3}}
\quad \mbox{ and } \quad
\widetilde  \psi_j'(t) = - \frac{\lambda^{5/3}\beta \widetilde \varphi_j}{(\widetilde  \varphi_j^2 + \widetilde  \psi_j^2)^{1/3}},
$$
where $\lambda$ is a large, positive constant defined later. One can check that $\widetilde  \varphi_j (t) = \lambda \varphi(\lambda t)$ and $\widetilde  \psi_j (t) = \lambda \psi(\lambda t)$ where $\varphi_j$ and $\psi_j$ are solutions of \eqref{eq-vpj} and
\begin{equation}\label{cd-vpj-m}
\varphi_j(0) + \psi_j(0) = \lambda^{-1} a,  \quad - \alpha \varphi_j(0) - \beta \psi_j (0) = 0,
\end{equation}
instead of \eqref{cd-vpj}. One then can obtain the dynamics for $\eta_j$ by choosing $\lambda$ large enough.\qed

\medskip
\noindent \textbf{Acknowledgments.} The authors were partially supported by  ANR Finite4SoS ANR-15-CE23-0007. H.-M. Nguyen thanks Fondation des Sciences Math\'ematiques de Paris (FSMP) for the Chaire d'excellence which allows him to visit  Laboratoire Jacques Louis Lions  and Mines ParisTech. This work has been done during this visit.

 \providecommand{\bysame}{\leavevmode\hbox to3em{\hrulefill}\thinspace}
\providecommand{\MR}{\relax\ifhmode\unskip\space\fi MR }
\providecommand{\MRhref}[2]{%
  \href{http://www.ams.org/mathscinet-getitem?mr=#1}{#2}
}
\providecommand{\href}[2]{#2}

\end{document}